
\input amssym  

\newdimen\normalparindent

\newdimen\paperwidth
\newdimen\paperheight
\paperwidth=210truemm
\paperheight=297truemm


\iffalse  

\hsize=15truecm
\hoffset=.46truecm
\vsize=23.7truecm
\voffset=.46truecm

\normalparindent=24pt

\font\elevenrm=cmr10 at 11pt 
\font\eightrm=cmr8 
\font\sixrm=cmr6 

\font\eleveni=cmmi10 at 11pt 
\font\eighti=cmmi8
\font\sixi=cmmi6

\font\elevensy=cmsy10 at 11pt 
\font\eightsy=cmsy8
\font\sixsy=cmsy6

\font\elevenex=cmex10 at 11pt 

\font\elevenbf=cmbx10 at 11pt 
\font\eightbf=cmbx8
\font\sixbf=cmbx6

\font\eleventt=cmtt10 at 11pt 
\font\elevensl=cmsl10 at 11pt 
\font\elevenit=cmti10 at 11pt 

\textfont0=\elevenrm \scriptfont0=\eightrm \scriptscriptfont0=\sixrm
\def\rm{\fam0\elevenrm}
\textfont1=\eleveni \scriptfont1=\eighti \scriptscriptfont1=\sixi
 
\textfont2=\elevensy \scriptfont2=\eightsy \scriptscriptfont2=\sixsy
\def\cal{\fam2}
\textfont3=\elevenex \scriptfont3=\elevenex \scriptscriptfont3=\elevenex
\textfont\itfam=\elevenit
\def\it{\fam\itfam\elevenit}
\textfont\slfam=\elevensl
\def\sl{\fam\slfam\elevensl}
\textfont\bffam=\elevenbf \scriptfont\bffam=\eightbf
\scriptscriptfont\bffam=\sixbf
\def\bf{\fam\bffam\elevenbf}
\textfont\ttfam=\eleventt
\def\tt{\fam\ttfam\eleventt}

\skewchar\eleveni='177 \skewchar\eighti='177 \skewchar\sixi='177
\skewchar\elevensy='60 \skewchar\eightsy='60 \skewchar\sixsy='60

\font\sc=cmcsc10 at 11pt 
\font\titlefont=cmssbx10 scaled \magstep3


\smallskipamount=3.5pt plus 1pt minus 1pt
\medskipamount=7pt plus 2pt minus 2pt
\bigskipamount=14pt plus 2pt minus 2pt
\normalbaselineskip=14pt
\normallineskip=1pt
\normallineskiplimit=0pt
\jot=3.5pt

\normalbaselines
\rm

\else  

\hsize=15truecm
\hoffset=.46truecm
\vsize=23.7truecm
\voffset=.46truecm

\normalparindent=20pt

\font\sc=cmcsc10


\font\titlefont=cmssbx10 scaled \magstep2

\fi

\def\S{\mathhexbox278\thinspace}
\def\SS{\mathhexbox278\mathhexbox278\thinspace}

\def\square{\hbox to.77778em{%
\hfil\vrule\vbox to.675em{\hrule width.6em\vfil\hrule}\vrule\hfil}}

\def\definition#1\par{\medbreak\noindent{\bf Definition.}\enspace
  #1\par\medbreak}
\def\example#1\par{\medbreak\noindent{\bf Example.}\enspace
  #1\par\medbreak}
\long\def\remark#1\par{\medbreak\noindent{\it Remark\/}.\enspace#1\par\medbreak}
\def\exercise#1\par{\medbreak\noindent{\bf Exercise.}\enspace
#1\par\medbreak}
\def\notation#1\par{\medbreak\noindent{\bf Notation.}\enspace
#1\par\medbreak}
\def\proof{\noindent{\it Proof\/}.\enspace}
\def\endproof{\nobreak\hfill\quad\square\par\medbreak}

\def\lineover#1{{\offinterlineskip\mathchoice
{\setbox0=\hbox{$\displaystyle#1$}%
\vbox{\kern .33pt\hbox to\wd0{\kern 1pt\leaders\hrule height .33pt%
\hfill\kern 1pt}\kern 1pt\box0}}
{\setbox0=\hbox{$\textstyle#1$}%
\vbox{\kern .33pt\hbox to\wd0{\kern 1pt\leaders\hrule height .33pt%
\hfill\kern 1pt}\kern 1pt\box0}}
{\setbox0=\hbox{$\scriptstyle#1$}%
\vbox{\kern .25pt\hbox to\wd0{\kern .8pt\leaders\hrule height .25pt%
\hfill\kern .8pt}\kern .8pt\box0}}
{\setbox0=\hbox{$\scriptscriptstyle#1$}%
\vbox{\kern .2pt\hbox to\wd0{\kern .6pt\leaders\hrule height .2pt%
\hfill\kern .6pt}\kern .6pt\box0}}}}

\def\Aut{\mathop{\rm Aut}\nolimits}

\def\End{\mathop{\rm End}\nolimits}

\def\Gal{\mathop{\rm Gal}\nolimits}
\def\GL{\mathop{\rm GL}\nolimits}
\def\Hom{\mathop{\rm Hom}\nolimits}
\def\id{{\rm id}}

\def\Spec{\mathop{\rm Spec}\nolimits}

\def\isom{\buildrel\sim\over\longrightarrow}
\def\morphism#1{\buildrel#1\over\longrightarrow}
\def\isomorphism#1{\mathrel{\mathop{\longrightarrow}%
\limits^{#1}_{\raise0.5ex\hbox{$\scriptstyle\sim$}}}}

\def\injlim{\mathop{\vtop{\offinterlineskip\halign{##\cr
 \hfil\rm lim\hfil\cr\noalign{\kern.1ex}\rightarrowfill\cr
 \noalign{\kern-.4ex}\cr}}}}
\def\projlim{\mathop{\vtop{\offinterlineskip\halign{##\cr
 \hfil\rm lim\hfil\cr\noalign{\kern.1ex}\leftarrowfill\cr
 \noalign{\kern-.4ex}\cr}}}}

\def\blank{\mkern12mu}

\def\textfrac#1/#2{{\textstyle{#1\over#2}}}

\def\commdiag#1{{
\def\rightar{\longrightarrow}
\def\downar{\big\downarrow}
\def\upar{\big\uparrow}
\def\hookrightar{\lhook\joinrel\longrightarrow}
\def\rightlabel##1{\rlap{$\scriptstyle##1$}}
\def\leftlabel##1{\llap{$\scriptstyle##1$}}
\vcenter{\baselineskip=0pt \lineskiplimit=0pt \lineskip=6pt
\halign{\hfil$\displaystyle##$\hfil&&\enskip\hfil$\displaystyle##$\hfil\crcr
#1\crcr}}}}

\def\relativediag#1#2#3#4#5#6{\vcenter{\baselineskip=3ex \halign{
\hfil$##$&$##$&$##$\hfil\cr
#1\quad& \hfilneg\buildrel#2\over\longrightarrow\hfilneg& \quad#3\cr
\lower1ex\llap{$\scriptstyle#4\hskip-1ex$}\searrow& \quad&
\swarrow\lower1ex\rlap{$\hskip-1ex\scriptstyle#5$} \cr
& \hfilneg#6\hfilneg&\cr}}}

\def\trianglediag#1#2#3#4#5#6{\vcenter{\baselineskip=3ex \halign{
\hfil$##$&$##$&$##$\hfil\cr
#1\quad& \hfilneg\buildrel#2\over\longrightarrow\hfilneg& \quad#3\cr
\lower1ex\llap{$\scriptstyle#6\hskip-1ex$}\nwarrow& \quad&
\swarrow\lower1ex\rlap{$\hskip-1ex\scriptstyle#4$} \cr
& \hfilneg#5\hfilneg&\cr}}}

\def\correspondence#1#2#3#4#5{\vcenter{\baselineskip=3ex \halign{
\hfil$##$&$##$&$##$\hfil\cr
&\hfilneg#1\hfilneg\cr
\raise1ex\llap{$\scriptstyle#2$}\swarrow&&\searrow
\raise1ex\rlap{$\scriptstyle#3$}\cr
#4&&#5\cr}}}

\newif\iffirstpar
\everypar{\iffirstpar\parindent=\normalparindent\firstparfalse\fi}

\def\sectionheading#1{\subcount=0 \subsectioncount=0 \eqcount=0
  \bigskip\vskip\parskip
  \leftline{\bf #1}\nobreak\smallskip\firstpartrue\parindent=0pt}

\def\section#1\par{\advance\sectioncount by1%
  \edef\currentlabel{\number\sectioncount}%
  \sectionheading{\number\sectioncount.\enspace#1}}

\def\unnumberedsection#1\par{\sectionheading{#1}}

\def\subsection#1\par{\medbreak\penalty-200\advance\subsectioncount by1%
  \edef\currentlabel{\number\sectioncount.\number\subsectioncount}%
  \leftline{\it\number\sectioncount.\number\subsectioncount.\enspace#1}%
  \smallskip\parindent=0pt\firstpartrue}

\newwrite\auxfile

\newcount\sectioncount \sectioncount=0
\newcount\subsectioncount 
\newcount\subcount 
\newcount\eqcount 

\def\subno{\global\advance\subcount by1\relax
  \number\sectioncount.\number\subcount
  \xdef\currentlabel{\number\sectioncount.\number\subcount}}
\def\proclaim #1. #2\par{\medbreak
  \noindent{\bf#1~\subno.\enspace}{\sl#2\par}%
  \ifdim\lastskip<\medskipamount \removelastskip\penalty55\medskip\fi}
\def\proclaimx #1 (#2). #3\par{\medbreak
  \noindent{\bf#1~\subno\ \rm (#2).\enspace}{\sl#3\par}%
  \ifdim\lastskip<\medskipamount \removelastskip\penalty55\medskip\fi}

\newdimen\algindent
\def\plusindent{\advance\algindent by \parindent}
\def\minusindent{\advance\algindent by-\parindent}


\newcount\algstepcount

\long\def\algorithm (#1). #2\endalgorithm{\medbreak
  \algindent=0pt%
  \algstepcount=0%
  \noindent{\bf Algorithm~\subno} (#1). {\sl#2}\par\medbreak}

\def\step{\advance\algstepcount by1
\edef\currentlabel{\number\algstepcount}
\smallskip\hangindent\parindent
\advance\hangindent by\algindent\indent
\llap{{\bf \the\algstepcount.}\enspace}\kern\algindent
\ignorespaces}


\def\labeldef#1#2{\expandafter\gdef\csname L@#1\endcsname{#2}}
\def\label#1{%
  \expandafter\xdef\csname L@#1\endcsname{\currentlabel}%
  \write\auxfile{\string\labeldef{#1}{\csname L@#1\endcsname}}%
  \ignorespaces}
\def\ref#1{\expandafter\ifx\csname L@#1\endcsname\relax
  \message{Undefined label `#1'}??\else
  \csname L@#1\endcsname\fi}

\def\eqnumber#1{\global\advance\eqcount by1\relax
  \eqno(\number\sectioncount.\number\eqcount)%
  \expandafter\xdef\csname E@#1\endcsname{%
    \number\sectioncount.\number\eqcount}}
\def\eqref#1{\expandafter\ifx\csname E@#1\endcsname\relax
  \message{Undefined equation `#1'}??\else
  (\csname E@#1\endcsname)\fi}

\newcount\refcount \refcount=0
\def\citedef#1#2{\expandafter\gdef\csname C@#1\endcsname{#2}}
\def\cite#1{\expandafter\ifx\csname C@#1\endcsname\relax
  \message{Undefined reference `#1'}\citedef{#1}{??}\fi
  \expandafter\gdef\csname R@#1\endcsname{\relax}%
  [\csname C@#1\endcsname]}
\def\citex#1#2{\expandafter\ifx\csname C@#1\endcsname\relax
  \message{Undefined reference `#1'}\citedef{#1}{??}\fi
  \expandafter\gdef\csname R@#1\endcsname{\relax}%
  [\csname C@#1\endcsname, #2]}
\def\reference#1{\advance\refcount by 1%
  \expandafter\ifx\csname R@#1\endcsname\relax
  \message{Warning: reference `#1' not used}\fi
  \expandafter\edef\csname C@#1\endcsname{\the\refcount}%
  \write\auxfile{\string\citedef{#1}{\csname C@#1\endcsname}}%
  \item{[\csname C@#1\endcsname]}}

\newif\ifauxexists
\immediate\openin0=\jobname.aux
\ifeof 0
  \auxexistsfalse
\else
  \auxexiststrue
\fi
\immediate\closein0
\ifauxexists
  \input \jobname.aux
\else
  \message{No file `\jobname.aux'}
\fi
\openout\auxfile=\jobname.aux

\def\Alg{\mathop{\bf Alg}\nolimits}
\def\C{{\bf C}}
\def\D{{\rm D}}
\def\DP{\mathop{\bf DP}\nolimits}
\def\F{{\bf F}}
\def\G{{\rm G}}
\def\Gm{{\bf G}_{\rm m}}
\def\Gmover#1{{\bf G}_{{\rm m},#1}}
\def\GS{\mathop{\bf GS}\nolimits}
\def\H{{\rm H}}
\def\HA{\mathop{\bf HA}\nolimits}
\def\sHom{\mathop{\bf Hom}\nolimits}
\def\Mod{\mathop{\bf Mod}\nolimits}
\def\O{{\cal O}}
\def\fp{{\frak p}}
\def\Q{{\bf Q}}
\def\Qbar{\overline{\bf Q}}
\def\R{{\bf R}}
\def\X{{\rm X}}
\def\Z{{\bf Z}}

\def\al{{\rm al}}
\def\op{{\rm op}}

\let\pairing\innerprod

\centerline{\titlefont Dual pairs of algebras and
finite commutative group schemes}

\bigskip

\centerline{Peter Bruin\footnote{}{%
The author was supported by a Veni grant from the Netherlands
Organisation for Scientific Research (NWO).}}
\smallskip
\centerline{28 September 2017}

\bigskip

{\narrower\narrower\noindent{\it Abstract.\/} We introduce a category
of dual pairs of finite locally free algebras over a ring.  This gives
an efficient way to represent finite locally free commutative group
schemes.  We give a number of algorithms to compute with dual pairs of
algebras, and we apply our results to Galois representations on finite
Abelian groups.\par}

\section Introduction

In this article we introduce an efficient way to write down finite
locally free commutative group schemes over a ring~$R$.  Our approach
is motivated by Cartier duality and is based on simultaneously
representing a group scheme and its Cartier dual.  This gives rise to
the definition of a category of {\it dual pairs of finite locally free
$R$-algebras\/}.  On the theoretical side, we prove that this category
is equivalent to the category of finite locally free commutative group
schemes over~$R$, and anti-equivalent to the corresponding category of
finite locally free commutative cocommutative Hopf algebras over~$R$.
On the algorithmic side, we show that this way of representing group
schemes leads to transparent algorithms for performing various
important operations with group schemes and groups of $S$-valued
points for $R$-algebras $S$.

In the case $R=\Q$, dual pairs of $\Q$-algebras give a concise way to
write down representations of the absolute Galois group of~$\Q$ on
finite Abelian groups.  An important advantage of our approach is that
in practice, the data that one needs to compute with has considerably
smaller height than when the computations are done using Hopf
algebras.  This is especially relevant when computing the data
numerically, such as in the explicit computations of Galois
representations by Bosman (see \cite{Bosman}
and \citex{Edixhoven-Couveignes}{Chapters 6 and~7}), Tian \cite{Tian},
Yin and Zeng \cite{Yin-Zeng}, Derickx, van Hoeij and Zeng \cite{DvHZ},
Mascot \cite{Mascot}, and the author (unpublished, see \cite{thesis}
and~\cite{modgalrep}).  Until now, the output of these computations
has consisted of polynomials whose splitting fields cut out the
sought-for Galois representations.  When one has computed such a
polynomial $F$ but no further data, one has to certify that $F$ has
the expected Galois group and ramification properties in order to
prove that $F$ indeed cuts out the expected Galois representation.
This approach does yield enough information to compute conjugacy
classes of Frobenius elements, but does not easily produce data that
amounts to specifying the group scheme structure.  The approach
described here makes it easier to write down such data, namely an
object of a category that is equivalent to the category of finite
commutative group schemes over~$\Q$.  In particular, the fact that a
dual pair of algebras possesses more structure than just the defining
polynomials of the algebras can be exploited to make the certification
of the Galois group much more direct.

The outline of the paper is as follows.
In~\S\ref{sec:notation-conventions} we fix the general notation and
conventions that we will use.  In \S\ref{sec:dual-pairs}, we define
the category of dual pairs of algebras over a ring~$R$ and compare it
(via the category of Hopf algebras over~$R$) to the category of finite
locally free commutative group schemes over~$R$.
In \S\ref{sec:examples}, we give two examples (one \'etale and one
non-\'etale) of dual pairs of algebras representing the 2-torsion
subschemes of certain elliptic curves.  In \S\ref{sec:algorithms}, we
describe algorithms for a number of relevant operations with group
schemes and their groups of points in our setting.
In \S\ref{sec:galois-representations}, we focus on Galois
representations of a field~$K$, which we can view as finite
commutative group schemes over~$K$.  We relate our description of
Galois representations to other descriptions and give an example to
show that our description leads to data of particularly small height.
In \S\ref{sec:future}, we sketch two directions for future work.  In
the appendix, we give a generic algorithm to identify a finite Abelian
group from a certain type of ``pairing matrix'' that is used
in \S\ref{sec:algorithms}.

\remark For most of the article, we work with group schemes over a
base ring $R$.  Whenever we are in this setting, all our constructions
are compatible with arbitrary base change and can be extended without
difficulties to arbitrary base schemes.

\goodbreak

\section Notation and conventions

\label{sec:notation-conventions}

All rings (and in particular all algebras) are assumed to be
commutative, unless otherwise stated.  Similarly, all Hopf algebras
are assumed to be commutative and cocommutative.

If $M$ and~$N$ are modules over a ring~$R$, a bilinear map $M\times
N\to R$ is said to be {\it perfect\/} if the induced $R$-bilinear maps
$$
\eqalign{
M&\longrightarrow\Hom_{\Mod_R}(N,R),\cr
N&\longrightarrow\Hom_{\Mod_R}(M,R)}
$$
are isomorphisms.

Let $R$ be a ring, and let $M$ be a finite locally free $R$-module.
We write $M^\vee=\Hom_{\Mod_R}(M,R)$ for the $R$-linear dual of~$M$.
We have a canonical perfect $R$-bilinear map
$$
\Phi_M\colon M\times M^\vee\to R.
$$
Furthermore, the canonical $R$-linear map
$$
\eqalign{
M\otimes_R M^\vee&\longrightarrow\Hom_{\Mod_R}(M,M)\cr
m\otimes\phi&\longmapsto(x\mapsto\phi(x)m)}
$$
is an isomorphism, and we write
$$
\theta_M \in M\otimes_R M^\vee
$$
for the unique element mapping to $\id_M\in\Hom_{\Mod_R}(M,M)$ under
the above isomorphism.  If $M$ is free over~$R$, then given an
$R$-basis $(a_1,\ldots,a_n)$ of~$M$ and the corresponding dual basis
$(b_1,\ldots,b_n)$ of $M^\vee$, we have
$$
\theta_M = \sum_{i=1}^n a_i\otimes b_i.
$$

Let $M$ and~$N$ be finite free modules over a ring~$R$ with bases
$(x_1,\ldots,x_m)$ and $(y_1,\ldots,y_n)$, respectively.  We will
sometimes identify $M\otimes_R N$ with the $R$-module of $m\times
n$-matrices with coefficients in~$R$ by representing an element
$\sum_{i=1}^m\sum_{j=1}^n c_{i,j} x_i\otimes y_j$ as the matrix
$(c_{i,j})_{1\le i\le m\atop 1\le j\le n}$.  Similarly, if $\Phi\colon
M\times N\to R$ is an $R$-bilinear map, we will sometimes represent
$\Phi$ by the $m\times n$-matrix $(\Phi(x_i,y_j))_{1\le i\le m\atop
1\le j\le n}$.  Furthermore, we will write $\Phi^{\rm t}$ for the
bilinear map $N\times M\to R$ defined by $(n,m)\mapsto\Phi(m,n)$.

Let $M$ and~$N$ be finite locally free modules over a ring~$R$, and
let $\Phi\colon M\times N\to R$ be a perfect $R$-bilinear map.  By
assumption, the $R$-linear map
$$
\eqalign{
N&\isom M^\vee\cr
n&\longmapsto(m\mapsto\Phi(m,n))}
$$
is an isomorphism.  Inverting this, we obtain an isomorphism
$M^\vee\isom N$, which in turn determines an element
$$
\theta_\Phi\in M\otimes_R N
$$
via the isomorphism
$$
\eqalign{
M\otimes_R N&\isom\Hom(M^\vee,N)\cr
m\otimes n&\longmapsto(\phi\mapsto\phi(m)n).}
$$
When $M$ and~$N$ are free $R$-modules with given bases, one checks
easily that the matrices of $\Phi$ and~$\theta_\Phi$ with respect to
these bases are inverse transposes of each other.

\goodbreak

\section Dual pairs of algebras

\label{sec:dual-pairs}

\subsection Motivation: Cartier duality for group schemes and Hopf
algebras

Let $H$ be a finite locally free commutative group scheme over a
ring~$R$, and let $A$ be the corresponding finite locally free Hopf
algebra over~$R$, so $A\cong\O_H(H)$ and $H\cong\Spec A$.  We denote
the structure maps of the Hopf algebra~$A$ by
$$
\displaylines{
m\colon A\otimes_R A\to A\quad\hbox{(multiplication)},\quad
e\colon R\to A\quad\hbox{(unit)},\cr
\mu\colon A\to A\otimes_R A\quad\hbox{(comultiplication)},\quad
\epsilon\colon A\to R\quad\hbox{(counit)}.}
$$

\remark Part of the definition of a Hopf algebra is the existence of
an antipode $A\isom A$.  However, the antipode itself does not need to
be included in the data defining the Hopf algebra structure, since an
antipode, if it exists, is unique.

Let $A^\vee$ be the Hopf algebra dual to~$A$; as an $R$-module, this
is defined by
$$
A^\vee = \Hom_{\Mod_R}(A,R)
$$
and the structure maps of~$A^\vee$ are defined by dualising those
of~$A$.  Let $H^*$ be the Cartier dual of~$H$, i.e.\ the spectrum
of~$A^\vee$; this is again a finite locally free commutative group
scheme over~$R$.  It is well known that $H^*$ represents the sheaf of
Abelian groups $\sHom(H,\Gmover R)$ for the {\it fppf\/} topology on
$\Spec R$; see for example Oort~\citex{Oort}{Theorem~16.1}.  We
therefore have a canonical morphism of $R$-schemes
$$
H\times H^*\longrightarrow\Gmover R
$$
and a corresponding $R$-algebra homomorphism
$$
R[x,x^{-1}]\longrightarrow A\otimes_R A^\vee.
$$
The element $\theta_A$ defined earlier is the image of~$x$ under this
homomorphism.  Furthermore, a theorem of Deligne stating that a finite
locally free commutative group scheme is annihilated by its rank (see
Tate and Oort \citex{Tate-Oort}{\S1}) implies that $\theta_A$ is a
(not necessarily primitive) $n$-th root of unity, where $n$ is the
rank of~$A$.

The basic observation that motivates our definition of dual pairs of
algebras in \S\ref{subsec:dual-pairs} below is that if $A$ is a finite
locally free Hopf algebra over~$R$, then $A$ is determined up to
isomorphism by the triple $(A^\al,(A^\vee)^\al,\Phi)$, where the
subscript ``al'' indicates the forgetful functor from the category of
Hopf algebras over~$R$ to the category of algebras over~$R$ and $\Phi$
is the canonical $R$-bilinear map $A\times A^\vee\to R$.

\subsection The category of dual pairs of algebras over a ring

\label{subsec:dual-pairs}

We introduce some notation in preparation for our main definition.  If
$A$ and~$B$ are two finite locally free $R$-algebras and $\Phi\colon
A\times B\to R$ is a perfect $R$-bilinear map, then $\Phi$ induces a
perfect $R$-bilinear map
$$
\eqalign{
\Phi^{(2)}\colon(A\otimes_R A)\times(B\otimes_R B)&\longrightarrow R\cr
(a\otimes a',b\otimes b')&\longmapsto\Phi(a,b)\Phi(a',b').}
$$
We define
$$
\mu_1^\Phi\colon A\to A\otimes_R A
\quad\hbox{and}\quad
\mu_2^\Phi\colon B\to B\otimes_R B
$$
as the unique $R$-linear maps making the diagrams
$$
\commdiag{
A & \isom & B^\vee\cr
\leftlabel{\mu_1^\Phi}\downar & & \downar\cr
A\otimes_R A & \smash\isom & (B\otimes_R B)^\vee}
\qquad\hbox{and}\qquad
\commdiag{
B & \isom & A^\vee\cr
\leftlabel{\mu_2^\Phi}\downar & & \downar\cr
B\otimes_R B & \smash\isom & (A\otimes_R A)^\vee}
$$
commutative, where the top horizontal arrows are induced by $\Phi$,
the bottom horizontal arrows are induced by $\Phi^{(2)}$ and the right
vertical arrows are the duals of the multiplication maps of $B$
and~$A$, respectively.  In other words, $\mu_1^\Phi$ and $\mu_2^\Phi$
are uniquely determined by the identities
$$
\eqalign{
\Phi(a,bb') &= \Phi^{(2)}(\mu_1^\Phi(a),b\otimes b')
\quad\hbox{for all $a\in A$, $b,b'\in B$},\cr
\Phi(aa',b) &= \Phi^{(2)}(a\otimes a',\mu_2^\Phi(b))
\quad\hbox{for all $a,a'\in A$, $b\in B$}.}
$$

Analogously, we define
$$
\epsilon_1^\Phi\colon A\to R
\quad\hbox{and}\quad
\epsilon_2^\Phi\colon B\to R
$$
as the unique $R$-linear maps making the diagrams
$$
\commdiag{
A & \isom & B^\vee\cr
\leftlabel{\epsilon_1^\Phi}\downar & & \downar\cr
R & \smash{\morphism{\id}} & R}
\qquad\hbox{and}\qquad
\commdiag{
B & \isom & A^\vee\cr
\leftlabel{\epsilon_2^\Phi}\downar & & \downar\cr
R & \smash{\morphism{\id}} & R}
$$
commutative, where the right vertical arrows are defined as evaluation
in the unit elements of $B$ and~$A$, respectively.  Thus
$\epsilon_1^\Phi$ and~$\epsilon_2^\Phi$ are uniquely determined by the
identities
$$
\eqalign{
\epsilon_1^\Phi(a)&=\Phi(a,1_B)
\quad\hbox{for all }a\in A,\cr
\epsilon_2^\Phi(b)&=\Phi(1_A,b)
\quad\hbox{for all }b\in B.}
$$

\definition Let $R$ be a ring.  A {\it dual pair of algebras\/}
over~$R$ is a triple $(A,B,\Phi)$, where $A$ and~$B$ are finite
locally free $R$-algebras and where
$$
\Phi\colon A\times B\to R
$$
is a perfect $R$-bilinear map such that the following hold:
\smallskip
\item{(1)} we have $\Phi(1_A,1_B)=1$;
\smallskip
\item{(2)} for all $a,a'\in A$ we have
$\Phi(aa',1_B)=\Phi(a,1_B)\Phi(a',1_B)$;
\smallskip
\item{(3)} for all $b,b'\in B$ we have
$\Phi(1_A,bb')=\Phi(1_A,b)\Phi(1_A,b')$;
\smallskip
\item{(4)} for all $a,a'\in A$ and $b,b'\in B$ we have
$$
\eqalign{
\Phi^{(2)}(\mu_1^\Phi(a)\mu_1^\Phi(a'),b\otimes b')
&=\Phi(aa',bb')\cr
&=\Phi^{(2)}(a\otimes a',\mu_2^\Phi(b)\mu_2^\Phi(b')).}
$$

\remark The conditions (1)--(4) are equivalent to saying that the
$R$-linear maps $\mu_1^\Phi$, $\mu_2^\Phi$, $\epsilon_1^\Phi$
and~$\epsilon_2^\Phi$ are in fact homomorphisms of $R$-algebras.

\definition Let $(A,B,\Phi)$, $(A',B',\Phi')$ be two dual pairs of
algebras over a ring $R$.  A {\it morphism\/} from $(A,B,\Phi)$ to
$(A',B',\Phi')$ is a pair of $R$-algebra homomorphisms $(f\colon A'\to
A,g\colon B\to B')$ satisfying
$$
\Phi(f(a'),b) = \Phi'(a',g(b))
\quad\hbox{for all $a'\in A'$ and $b\in B$}.
\eqnumber{eq:Phi-adjoint}
$$

Morphisms can be composed as follows: given a morphism $(f,g)$ from
$(A,B,\Phi)$ to $(A',B',\Phi')$ and a morphism $(f',g')$ from
$(A',B',\Phi')$ to $(A'',B'',\Phi'')$, we put
$$
(f',g')\circ(f,g)=(f\circ f',g'\circ g);
$$
one checks immediately that this is a morphism from $(A,B,\Phi)$ to
$(A'',B'',\Phi'')$.

For every ring $R$, we denote by $\DP_R$ the category of dual pairs of
$R$-algebras with morphisms as defined above.

\remark Alternatively, we could have declared a morphism $(f,g)$ as
above to be a ``morphism from $(A',B',\Phi')$ to $(A,B,\Phi)$''
instead of the other way around.  The reason for our chosen convention
is that it gives us an equivalence of categories, rather than an
anti-equivalence, in Corollary~\ref{cor:equivalence}.

\subsection (Anti-)equivalences with Hopf algebras and group schemes

Let $R$ be a ring.  We write $\HA_R$ for the category of finite,
locally free, commutative and cocommutative Hopf algebras over~$R$.

If $(A,B,\Phi)$ is a dual pair of $R$-algebras, we write
$$
\H_R(A,B,\Phi) = (A,m_A,e_A,\mu_1^\Phi,\epsilon_1^\Phi);
$$
see \S\ref{subsec:dual-pairs} for the definition of $\mu_1^\Phi$ and
$\epsilon_1^\Phi$.  The definition of dual pairs implies that
$\H_R(A,B,\Phi)$ is an object of $\HA_R$.  Furthermore, given a
morphism $(f,g)\colon(A,B,\Phi)\to(A',B',\Phi)$ of dual pairs of
$R$-algebras, the morphism $f\colon A'\to A$ is a homomorphism of Hopf
algebras in the opposite direction.  In this way we obtain a functor
$$
\H_R\colon\DP_R\longrightarrow\HA_R^\op.
$$

Conversely, if $(A,m,e,\mu,\epsilon)$ is a Hopf algebra over~$R$, we
write
$$
\D_R(A,m,e,\mu,\epsilon) = (A,A^\vee,\Phi_A).
$$
Here $A$ is equipped with the algebra structure defined by $m$
and~$e$, and $A^\vee$ is equipped with the algebra structure defined
by $\epsilon^\vee$ and $\mu^\vee$.  The definition of dual pairs
implies that $\D_R(A,m,e,\mu,\epsilon)$ is a dual pair of algebras
over~$R$.  Furthermore, given a homomorphism
$f\colon(A',m',e',\mu',\epsilon')\to(A,m,e,\mu,\epsilon)$ of Hopf
algebras, the pair $(f\colon A'\to A,f^\vee\colon A^\vee\to A'^\vee)$
is a morphism from the dual pair $(A,A^\vee,\Phi_A)$ to the dual pair
$(A',A'^\vee,\Phi_{A'})$.  In this way we obtain a functor
$$
\D_R\colon\HA_R^\op\longrightarrow\DP_R.
$$

\proclaim Theorem. Let $R$ be a ring.  The functors
$\H_R\colon\DP_R\to\HA_R^\op$ and $\D_R\colon\HA_R^\op\to\DP_R$ are
anti-equivalences of categories between the category of dual pairs
over~$R$ and the category of finite locally free commutative
cocommutative Hopf algebras over~$R$.

\label{theorem:equivalence}

\proof It follows immediately from the definitions that the functor
$\H_R\circ\D_R$ is the identity on the category $\HA_R^\op$.  The
functor $\D_R\circ\H_R$ sends a dual pair $(A,B,\Phi)$ to the dual
pair $(A,A^\vee,\Phi_A)$; there is a natural isomorphism
$$
\eta_{(A,B,\Phi)}\colon(A,B,\Phi)\isom(A,A^\vee,\Phi_A)
$$
in $\DP_R$ defined by the pair of isomorphisms $(\id\colon A\to
A,\phi\colon B\to A^\vee)$, where $\phi$ is defined by
$\phi(b)(a)=\Phi(a,b)$.  The isomorphisms $\eta_{A,B,\Phi}$ define an
isomorphism
$$
\eta\colon\id_{\DP_R}\isom\D_R\circ\H_R
$$
of functors from the category $\DP_R$ to itself.  One checks easily
that each of the natural isomorphisms
$$
\eqalign{
\eta\D_R\colon\D_R&\longrightarrow\D_R\circ\H_R\circ\D_R,\cr
\H_R\eta\colon\H_R&\longrightarrow\H_R\circ\D_R\circ\H_R}
$$
is the identity.  We conclude that $(\H_R,\D_R,\eta,\id)$ is an
adjoint equivalence of categories from $\DP_R$ to $\HA_R^\op$ (see
Mac~Lane \citex{Mac Lane}{\S IV.4} for the definition of an adjoint
equivalence).\endproof

\remark One may wonder why we have introduced the notion of dual pairs
when Theorem~\ref{theorem:equivalence} shows that it is essentially
the same as that of Hopf algebras.  The reason is that for algorithmic
purposes, we would like to avoid computing directly with
comultiplication maps.  The proof of Theorem~\ref{theorem:equivalence}
shows that $\HA_R^\op$ is canonically embedded into $\DP_R$, but the
category $\DP_R$ has more objects because in a dual pair $(A,B,\Phi)$
the $R$-module $B$ is only isomorphic, and not necessarily identical,
to~$A^\vee$.  We can use this extra ``degree of freedom'' to present
$B$ in a computationally efficient way; the multiplication map on~$B$
is then used to construct the comultiplication map on~$A$, rather than
the other way around.

\goodbreak

Let $\GS_R$ denote the category of finite locally free commutative
group schemes over~$R$.

\proclaim Corollary. Let $R$ be a ring.  There is a canonical
equivalence of categories
$$
\G_R\colon\DP_R\to\GS_R
$$
such that for every dual pair $(A,B,\Phi)$ of $R$-algebras, the
underlying $R$-scheme of\/ $\G_R(A,B,\Phi)$ equals $\Spec A$.

\label{cor:equivalence}

By construction, the above equivalence is compatible with duality in
the following sense: if $(A,B,\Phi)$ is a dual pair of $R$-algebras
and $H$ is the group scheme determined by $(A,B,\Phi)$, then the
Cartier dual $H^*$ is canonically isomorphic to the group scheme
determined by $(B,A,\Phi^{\rm t})$.  The resulting isomorphisms
$B\isom A^\vee$ and $A\isom B^\vee$ of $R$-modules equal those arising
from~$\Phi$, and the image of~$x$ under the resulting $R$-algebra
homomorphism
$$
R[x,x^{-1}]\longrightarrow A\otimes_R B
$$
equals the element $\theta\in A\otimes_R B$ defined earlier.

\goodbreak

\section Examples

\label{sec:examples}

\subsection The 2-torsion of an elliptic curve over~$\Q$

\label{sec:ellQ}

We consider the elliptic curve $E$ over~$\Q$ given by a Weierstra\ss\
equation of the form
$$
E\colon y^2 = x^3 - ax
$$
with $a\in\Q^\times$.

\proclaim Proposition. The group scheme $E[2]$ over~$\Q$ can be
represented by the dual pair $(A,B,\Phi)$ of\/ $\Q$-algebras, where
$$
A=B=\Q\times\Q\times\Q[t]/(t^2-a)
$$
and the matrix of the $\Q$-bilinear map $\Phi$ with respect to the
basis $((1,0,0),(0,1,0),(0,0,1),(0,0,t))$ is
$$
\Phi = \pmatrix{1/4& 1/4& 1/2& 0\cr 1/4& 1/4& -1/2& 0\cr
1/2& -1/2& 0& 0\cr 0& 0& 0& a}.
$$

\proof The coordinate ring of $E[2]$ is $\Q\times\Spec\Q[x]/(x^3-ax)$,
where the first factor comes from the origin of~$E$ and the second
factor from the points of order~$2$.  Using the isomorphism
$$
\eqalign{
\Q[x]/(x^3-ax)&\isom\Q\times\Q[t]/(t^2-a)\cr
x&\longmapsto(0,t),}
$$
we identify this coordinate ring with $A$.  The group scheme $E[2]$ is
identified with its own Cartier dual by the Weil pairing $E[2]\times
E[2]\to\mu_{2,\Q}$.

The splitting field of $E[2]$ is $L=\Q(\sqrt{a})$.  The four points of
$E[2]$ over this field correspond to the four $\Q$-algebra
homomorphisms
$$
\eqalign{
p_0,p_1,p_2,p_3\colon A&\longrightarrow L\cr
(a,b,c+dt)&\longmapsto
\cases{a& if $i=0$,\cr
b& if $i=1$,\cr
c+d\sqrt{a}& if $i=2$,\cr
c-d\sqrt{a}& if $i=3$.}}
$$
The matrix of the induced $L$-algebra isomorphism
$(p_0,p_1,p_2,p_3)\colon A\otimes_\Q L\to L^4$ with respect to the
above basis of~$A$ and the standard basis of~$L^4$ is
$$
P = \pmatrix{
1& 0& 0& 0\cr
0& 1& 0& 0\cr
0& 0& 1& \sqrt{a}\cr
0& 0& 1& -\sqrt{a}}.
$$
For $i,j\in\{0,1,2,3\}$, we consider the $\Q$-algebra homomorphism
$$
\eqalign{
p_i\otimes p_j\colon A\otimes_\Q B&\longrightarrow L\cr
a\otimes b&\longmapsto p_i(a)p_j(b)}
$$
The Weil pairing is determined by a the unique element $\theta\in
A\otimes_\Q B$ such that
$$
(p_i\otimes p_j)(\theta) = \cases{
1& if $i=0$, $j=0$ or $i=j$\cr
-1& otherwise.}
$$
Writing $\Theta$ for the $4\times 4$-matrix of coefficients
of~$\theta$ with respect to the fixed basis of $A\otimes_\Q B$, we see
that the above condition is equivalent to the matrix equation
$$
P\Theta P^{\rm t} = \pmatrix{
1& 1& 1& 1\cr
1& 1& -1& -1\cr
1& -1& 1& -1\cr
1& -1& -1& 1}.
$$
The unique solution is
$$
\Theta = \pmatrix{
1& 1& 1& 0\cr
1& 1& -1& 0\cr
1& -1& 0& 0\cr
0& 0& 0& 1/a}.
$$
By our earlier remark that $\Phi=\Theta^{{\rm t},-1}$, the proposition
is proved.\endproof

\subsection The 2-torsion of a supersingular elliptic curve in
characteristic 2

We now consider a non-\'etale group scheme, namely the 2-torsion of
the supersingular elliptic curve $E$ over~$\F_2$ given by the
(projective) Weierstra\ss\ equation
$$
E\colon y^2z+yz^2=x^3.
$$
As an effective Cartier divisor, $E[2]$ equals the point $(0:1:0)$
with multiplicity~$4$.

\proclaim Proposition. The group scheme $E[2]$ over~$\F_2$ can be
represented by the dual pair $(A,B,\Phi)$ of $\F_2$-algebras, where
$$
A=B=\F_2[t]/(t^4)
$$
and the matrix of the $\F_2$-bilinear map $\Phi$ with respect to the
power basis $(1,t,t^2,t^3)$ is
$$
\Phi = \pmatrix{1& 0& 0& 0\cr 0& 0& 1& 0\cr 0& 1& 0& 0\cr 0& 0& 0& 1}.
$$

\proof A local parameter of~$E$ at~$(0:1:0)$ is $t=x/y$; this
identifies $E[2]$ as a scheme with $\Spec A$.  Let
$E_A=E\times_{\Spec\F_2}\Spec A$.  For any section $P\colon A\to E_A$,
the image of~$P$ is an relative effective Cartier divisor on~$E_A$,
which we denote by~$[P]$.  The ``universal 2-torsion point on~$E$'' is
the section
$$
P\colon\Spec A\to E_A
$$
that can be described in the affine patch $\{y=1\}$ as $P=(t:1:t^3)$.
Let $D$ be the relative Cartier divisor $[P]-[O]$ on~$E_A$.  Then $2D$
is a principal divisor.  Below we compute a relative Cartier divisor
$D'$ on~$E_A$ that is linearly equivalent to~$D$ but has disjoint
support from~$D$, as well as rational functions $f$ and~$f'$ on~$E_A$
(global sections of the sheaf of total quotient rings; see
Hartshorne \citex{Hartshorne}{\S II.6}) with divisors $D$ and~$D'$,
respectively.  We then pull back $D$ and~$f$ to $E_{A\otimes_{\F_2}
B}=E\times_{\Spec\F_2}(\Spec A\times_{\Spec\F_2}\Spec B)$ via the
first projection $\Spec A\times_{\Spec\F_2}\Spec B\to\Spec A$, and we
pull back $D'$ and~$f'$ via the second projection $\Spec
A\times_{\Spec\F_2}\Spec B\to\Spec B$, identifying $\Spec B$ with
$\Spec A$.  Then the Weil pairing is given by the element
$$
\theta = f'(D)/f(D')\in\mu_2(A\otimes_{\F_2}B).
$$

We will use the following sections $\Spec A\to E_A$:
$$
O=(0:1:0),\quad
P=(t:1:t^3),\quad
Q=(0:0:1),\quad
R=(t:t^3:1).
$$
The lines $x=ty$ and $x=tz$ intersect $E_A$ in the relative effective
Cartier divisors $[P]+[Q]+[(t:1+t^3:1)]$ and $[O]+[R]+[(t:1+t^3:1)]$,
respectively; this shows that $P+Q=R$ and that $D=[P]-[O]$ is linearly
equivalent to $D'=[R]-[Q]$.  Furthermore, $2D=2[P]-2[O]$ is the
divisor of the rational function $f=(z-t^2x)/z$, and $2D'=2[R]-2[Q]$
is the divisor of the rational function $f'=(y-t^2x)/y$.  Pulling the
above divisors and functions back to the base ring $A\otimes_{\F_2}B$
as described above and writing this ring as
$\F_2[t_1,t_2]/(t_1^4,t_2^4)$, we obtain divisors and functions
$$
\displaylines{
D=[(t_1:1:t_1^3)]-[(0:1:0)],\quad
D'=[(t_2:t_2^3:1)]-[(0:0:1)],\cr
f=(z-t_1^2x)/z,\quad
f'=(y-t_2^2x)/y.}
$$
Using these data, we compute the Weil pairing on~$E[2]$ as
$$
\eqalign{
\theta &= {f'(D)\over f(D')}\cr
&={f'(t_1:1:t_1^3)/f'(0:1:0)\over f(t_2:t_2^3:1)/f(0:0:1)}\cr
&={(1-t_2^2 t_1)/1\over(1-t_1^2t_2)/1}\cr
&=1+t_1t_2^2+t_1^2t_2+t_1^3t_2^3.}
$$
This means that the matrix $\Theta$ of coefficients of~$\theta$ with
respect to the power basis $(1,t,t^2,t^3)$ is
$$
\Theta = \pmatrix{1& 0& 0& 0\cr 0& 0& 1& 0\cr 0& 1& 0& 0\cr 0& 0& 0& 1}.
$$
By our earlier remark that $\Phi=\Theta^{{\rm t},-1}$, the proposition
is proved.\endproof

\remark The above result allows us to recover the comultiplication map
by computing the group operation as in \S\ref{subsec:group-operation}
for the ``universal pair of points'' with values in
$A\otimes_{\F_2}A\simeq\F_2[t_1,t_2]/(t_1^4,t_2^4)$.  The result is
the $\F_2$-algebra homomorphism
$$
\eqalign{
\mu_{E[2]}\colon \F_2[t]/(t^4)&\longrightarrow\F_2[t_1,t_2]/(t_1^4,t_2^4)\cr
t&\longmapsto t_1+t_2+t_1^2t_2^2.}
$$
As expected, this agrees with the formal group of~$E$ (modulo $t^4$)
as computed from the Weierstra\ss\ equation; see for example
Tate \citex{Tate}{\S3}.

\section Algorithms

\label{sec:algorithms}

In this section we describe a number of algorithms for working with
finite locally free commutative group schemes over a ring~$R$.  We
will represent finite group schemes over~$R$ as dual pairs of
$R$-algebras, and describe our algorithms in terms of these.  We also
analyse the running time of our algorithms in terms of the required
number of operations in~$R$, and where applicable also the complexity
of factoring univariate polynomials over~$R$.

Throughout this section, the base ring~$R$ will be assumed to be
equipped with an ``algorithmic representation'', in the sense that we
have a way to write down elements of~$R$ (for example as finite bit
strings), and that we can perform ring operations and equality testing
using this representation.  Whenever we refer to an $R$-algebra $S$,
we assume $S$ to be represented algorithmically in a similar way
as~$R$.

\subsection Efficiently representing group schemes when the algebras
are monogenic

Let $(A,B,\Phi)$ be a dual pair of $R$-algebras.  We make the
assumption that the algebras $A$ and~$B$ are monogenic over~$R$, and
that we have an explicit presentation
$$
A = R[x]/(f),\quad B = R[y]/(g)
$$
where $f\in R[x]$ and $g\in R[y]$ are monic polynomials, say of
degree~$n\ge0$.  This implies that $A$ and~$B$ are free over~$R$; we
fix $R$-bases $(a_0,\ldots,a_{n-1})$ and $(b_0,\ldots,b_{n-1})$ of $A$
and~$B$, respectively, by
$$
a_i = x^i \bmod f
\quad\hbox{and}\quad
b_i = y^i \bmod g.
$$
We represent the perfect bilinear map $\Phi$ by its matrix with
respect to these basis.  In this way we represent $A$ and~$B$ by $n$
elements of~$R$ each, and $\Phi$ by $n^2$ elements of~$R$.

\remark In the important case where $R=\Q$, the algebras $A$ and~$B$
are products of number fields.  To minimise the height of the entries
of the matrix of the bilinear map~$\Phi$, it is useful to construct
our chosen $\Q$-bases of $A$ and~$B$ from integral bases of these
number fields instead of bases consisting of powers of fixed primitive
elements.

Below we will sometimes need to extend the base ring of our dual
pairs.  If $S$ is an $R$-algebra and $A=R[x]/(f)$, then we represent
the tensor product $A\otimes_R S$ as
$$
A\otimes_R S = S[x]/(f),
$$
where $f$ is now viewed as an element of $S[x]$, and similarly
for~$B$.

Given monic polynomials $f\in R[x]$ and~$g\in R[y]$ of degree~$n$ and
an $n\times n$-matrix $M$, it is in principle straightforward to check
whether these data define a dual pair of $R$-algebras by verifying the
conditions (1)--(4) in the definition.  However, the required number
of operations in~$R$ quickly becomes large; we will explain
in \SS\ref{subsec:group-structure} and~\ref{subsec:complex} below how
this verification can be done faster in cases where $R$ is a domain
and we can compute or approximate a common splitting field for $A$
and~$B$.

\subsection Groups of points

\label{subsec:group-operation}

Let $G$ be a finite locally free commutative group scheme over a
ring~$R$, represented by a dual pair of $R$-algebras $(A,B,\Phi)$.
Let $S$ be an $R$-algebra.  We will denote the group operation on
$G(S)=\Hom_{\Spec R}(S,G)$ by $(p,q)\mapsto p*q$.

Using the maps
$$
G(S)=\Hom_{\Spec R}(S,G)\cong\Hom_{\Alg_R}(A,S)\subseteq
\Hom_{\Mod_R}(A,S),
$$
we can represent elements of the Abelian group $G(S)$ by $R$-linear
maps $A\to S$.  The group operation on~$G(S)$ can be computed as
follows.  Given two elements of~$G(S)$, represented as $R$-linear maps
$$
p,q\colon A\to S,
$$
we use $\Phi$ to convert $p$ and~$q$ to elements $\hat p$ and~$\hat q$
of the $R$-algebra $B\otimes_R S$.  We then compute the product $\hat
p\hat q\in B\otimes_R S$ and convert the result back to an $R$-algebra
homomorphism $p*q\colon A\to S$ corresponding to the product of $p$
and~$q$ under the group operation on $G(S)$.

Let $S$ be an $R$-algebra.  Suppose that given a monic polynomial
$f\in R[x]$ we can compute the set of roots of~$f$ in~$S$.  As noted
above, the group of points $G(S)$ is canonically isomorphic to
$\Hom_{\Alg_R}(A,S)$.  If $A=R[x]/(f)$, then determining $G(S)$ as a
set is therefore equivalent to computing the set of roots of $f$
in~$S$.  In \S\ref{subsec:group-structure} below, we will describe how
the group structure on $G(S)$ can be computed in the case where $K$ is
a field of characteristic not dividing the degree $n$ of~$G$ and $S$
is a finite extension of~$K$ containing an $n$-th root of unity.

\subsection Cartier duality

\label{subsec:pairing}

Let $G$ be a finite locally free commutative group scheme over a
ring~$R$, represented by a dual pair of $R$-algebras $(A,B,\Phi)$, and
let $S$ be an $R$-algebra.  Cartier duality gives an isomorphism
$$
\tau\colon G^*(S)\isom\Hom(G_S,\Gm)
$$
of Abelian groups.  Given an element $q\in G^*(S)$, viewed as an
$R$-linear map $q\colon B\to S$, we apply the $R$-linear map
$$
\id_A\otimes q\colon A\otimes_R B\longrightarrow A\otimes_R S
$$
to the element $\theta_\Phi\in A\otimes_R B$.  The result is the
invertible element in $A\otimes_R S$ corresponding to $\tau(q)$.

We can compute the duality pairing
$$
\pairing\colon G(S)\times G^*(S)\longrightarrow\Gm(S)=S^\times
$$
in a similar way.  Given elements $p\in G(S)$ and $q\in G^*(S)$,
represented as $R$-linear maps $p\colon A\to S$ and $q\colon B\to S$,
we apply the $R$-linear map
$$
p\otimes q\colon A\otimes_R B\longrightarrow S\otimes_R S
$$
to the element $\theta_\Phi\in A\otimes_R B$ and apply the
multiplication map $S\otimes_R S$ to the result to obtain $\langle
p,q\rangle$.  If $A$ and~$B$ are free over~$R$ and we represent
$\theta_\Phi$ by a matrix over~$R$ and $p,q$ as row vectors over~$S$,
then $\langle p,q\rangle$ is nothing but the product $p\theta_\Phi
q^{\rm t}\in S$.

\subsection Groups of morphisms in $\DP_R$

Let $(A,B,\Phi)$ and $(A',B',\Phi')$ be two dual pairs of algebras
over a ring~$R$.  We will now make the group operation in
$\Hom_{\DP_R}((A,B,\Phi),(A',B',\Phi'))$ explicit.  Given two
morphisms
$$
(f,g),(f',g')\colon(A,B,\Phi)\longrightarrow(A',B',\Phi'),
$$
we denote the sum $(f,g)+(f',g')$ by $(f'',g'')$.  We compute
$(f'',g'')$ as follows.  Writing $G'$ for the group scheme
$\G_R(A',B',\Phi')$, we view $f,f'\colon A'\to A$ as $A$-valued points
of $G'$.  We add these as in \S\ref{subsec:group-operation} to obtain
$f''$.  We then compute $g''\colon B\to B'$ as the unique $R$-linear
map satisfying the adjointness property~\eqref{eq:Phi-adjoint}; this
is automatically an $R$-algebra homomorphism.

Moreover, if $R$ and~$A$ are such that given a monic polynomial $h\in
R[x]$ we can compute the set of roots of~$f$ in~$A$, then we can
compute the set $\Hom_{\DP_R}((A,B,\Phi),(A',B',\Phi'))$ as follows.
We first compute the set $\Hom_{\Alg_R}(A',A)$ of $R$-algebra
homomorphisms $f\colon A'\to A$.  For each $f\in\Hom_{\Alg_R}(A',A)$,
we then compute the unique $R$-linear map $g\colon B\to B'$ satisfying
the adjointness property~\eqref{eq:Phi-adjoint}.  The pair $(f,g)$ is
in $\Hom_{\DP_R}((A,B,\Phi),(A',B',\Phi'))$ if and only if $g$ is also
an $R$-algebra homomorphism.

Finally, we note that computing
$\Hom_{\DP_R}((A,B,\Phi),(A',B',\Phi'))$ allows us in particular to
determine whether the two objects $(A,B,\Phi)$ and $(A',B',\Phi')$ of
$\DP_R$ are isomorphic.

\subsection Direct sums, kernels and cokernels

The category $\GS_R$ of finite locally free commutative group schemes
over a ring is an Abelian category in some cases, for example if $R$
is a field.  However, $\GS_R$ fails to be an Abelian category in
general, because quotients of finite locally free commutative group
schemes are not necessarily locally free.  For example, the unique
non-trivial homomorphism $\Z/2\Z\to\mu_2$ of group schemes over~$\Z$
has trivial kernel and cokernel in~$\GS_\Z$, yet it is not an
isomorphism.  However, $\GS_R$ is always an exact category in the
sense of Quillen \citex{Quillen}{\S2}, and in particular is an
additive category.

Let $(A,B,\Phi)$ and $(A',B',\Phi')$ be two dual pairs of algebras
over a ring~$R$.  We define a perfect $R$-bilinear map between
$A\otimes_R A'$ and $B\otimes_R B'$ by
$$
\eqalign{
\Phi''\colon(A\otimes_R A')\times(B\otimes_R B')\longrightarrow R\cr
(a\otimes a',b\otimes b')\longmapsto\Phi(a,b)\Phi(a',b').}
$$
Writing $A''=A\otimes_R A$ and $B''=B\otimes_R B'$, one can check that
$(A'',B'',\Phi'')$ is again a dual pair.  It follows from the
compatibility of Cartier duality with direct sums that
$(A'',B'',\Phi'')$, together with the ``obvious'' morphisms to and
from the dual pairs $(A,B,\Phi)$ and $(A',B',\Phi')$, is a direct sum
(and product) of these two dual pairs.

If $R$ is a field, then we can compute kernels and cokernels of
morphisms in $\DP_R$.  Given a morphism $(f,g)$ from a dual pair
$(A,B,\Phi)$ to a dual pair $(A',B',\Phi')$, let
$$
(k,l)\colon(A'',B'',\Phi'')\longrightarrow(A,B,\Phi)
$$
be the kernel of~$(f,g)$.  Also, let
$$
(f_0,g_0)\colon(A,B,\Phi)\longrightarrow(A',B',\Phi')
$$
be the zero morphism of dual pairs.  We can then compute the quotient
algebra $k\colon A\to A''$ as the coequaliser of $f,f_0\colon A'\to
A$, and the subalgebra $l\colon B''\to B$ as the equaliser of
$g,g_0\colon B\to B'$.  Then $B''$ and $\ker k$ are each other's
orthogonal complements under~$\Phi$, and $\Phi''$ is the bilinear map
induced by~$\Phi$ on $A''\otimes_R B''$.  Computing the cokernel of
$(f,g)$ is entirely analogous.

\subsection Validating input

\label{subsec:validating-input}

Let $R$ be a domain, and let $K$ be the field of fractions of~$R$.
Suppose we are given a triple $(A,B,\Phi)$, where $A$ and~$B$ are two
finite locally free $R$-algebras of degree~$n$ such that the
characteristic of~$K$ does not divide $n$, and $\Phi\colon A\times
B\to R$ is a perfect $R$-bilinear map.  By a {\it common splitting
field\/} for $A$ and~$B$, we will mean an extension field $L$ of~$K$
such that the $L$-algebras $A\otimes_R L$ and $B\otimes_R L$ are both
isomorphic to~$L^n$.  Under the assumption that we can find such a
field~$L$, we will show how to decide whether is $(A,B,\Phi)$ a dual
pair of $R$-algebras.

We first note that if $L$ is any extension field of~$K$, then a
necessary and sufficient condition for $(A,B,\Phi)$ to be a dual pair
of $R$-algebras is that the base change $(A_L,B_L,\Phi_L)$ to~$L$ is a
dual pair of $L$-algebras.  The reason for this is that under our
assumptions on $(A,B,\Phi)$, the statement that $(A,B,\Phi)$ is a dual
pair of $R$-algebras is equivalent to a list of identities involving
maps between finite locally free $R$-modules.  These can be checked
after base extension to~$L$ because the modules are flat over~$R$ and
the homomorphism $R\to L$ is injective.  From now on, we will
therefore work over a field~$L$ as above.

We assume furthermore (enlarging $L$ if necessary) that $L$ is a
common splitting field of $A$ and~$B$, that $L$ contains a root of
unity of order~$n$, and that we can explicitly determine such a root
of unity $\zeta$ as well as all elements of $\Hom_{\Alg_L}(A,L)$ and
$\Hom_{\Alg_L}(B,L)$.  We have a function of sets
$$
\pairing\colon\Hom_{\Alg_L}(A,L)\times\Hom_{\Alg_L}(B,L)
\longrightarrow L
$$
defined as in \S\ref{subsec:pairing}.  If $(A,B,\Phi)$ is a dual pair
of $R$-algebras corresponding to a group scheme~$G$, then $\pairing$
can be identified with the duality pairing $G(L)\times G^*(L)\to
L^\times$, but the definition makes sense without assuming that
$(A,B,\Phi)$ is a dual pair.

We now observe (using the notation and terminology of the appendix)
that $(A_L,B_L,\Phi_L)$ is a dual pair of $L$-algebras if and only if
it is isomorphic to the dual pair $((H_d)_L,(H_d^*)_L,\Phi_d)$ for
some sequence of elementary divisors $d$.  Here $(H_d)_L$ and
$(H_d^*)_L$ are the constant group schemes over~$L$ associated with
$H_d$ and~$H_d^*=\Hom(H_d,\Q/\Z)$, and $\Phi_d$ is defined by viewing
$(H_d^*)_L$ as the Cartier dual of $(H_d)_L$ via the isomorphism
${1\over n}\Z/\Z\isom\langle\zeta\rangle$ sending $1/n$ to~$\zeta$.

\algorithm (Input validation). Given a ring~$R$, two finite locally
free algebras $A$ and~$B$ of degree~$n$ over~$R$, a perfect
$R$-bilinear map $\Phi\colon A\times B\to R$ and a common splitting
field $L$ of $A$ and~$B$ containing a root of unity of order~$n$, this
algorithm outputs ``True'' if $(A,B,\Phi)$ is a dual pair of
$R$-algebras, and ``False'' otherwise.

\step Fix bijections
$\alpha\colon\{1,\ldots,n\}\isom\Hom_{\Alg_R}(A,L)$ and
$\beta\colon\{1,\ldots,n\}\isom\Hom_{\Alg_R}(B,L)$.

\step Check whether all elements $\langle\alpha(i),\beta(j)\rangle\in
L$ are $n$-th roots of unity; if not, output ``False'' and stop.

\step Using the group isomorphism
$$
\lambda\colon\langle\zeta\rangle\isom{1\over n}\Z/\Z
$$
sending $\zeta$ to $1/n$, compute the matrix
$$
T=\bigl(\lambda(\langle\alpha(i),\beta(j)\rangle)\bigr)_{i,j=1}^n
$$
with entries in ${1\over n}\Z/\Z$.  (Note that $\lambda$ can be
computed by listing all powers of~$\zeta$, for example.)

\step Using Algorithm~\ref{algorithm:group-structure}, check whether
$T$ represents an Abelian group.  If so, output ``True'', otherwise
output ``False''.

\endalgorithm

The above procedure for validating the input relies on explicitly
computing a splitting field for $(A,B,\Phi)$.  In practice, it is
desirable to perform this task {\it without\/} explicitly computing a
splitting field.  In the case of dual pairs of $\Q$-algebras, we will
give an alternative approach in \S\ref{subsec:complex} below.

\remark The above procedure remains valid when the assumption ``$L$
contains a root of unity of order~$n$'' is weakened to ``$L$ contains
a root of unity of order equal to the exponent of $G(\bar L)$'', where
$\bar L$ is a separable closure of~$L$.  This exponent may be smaller
than~$n$, but is in general not known in advance.

\subsection Determining group structures

\label{subsec:group-structure}

Let $K$ be a field, and let $(A,B,\Phi)$ be a dual pair of
$K$-algebras corresponding to a finite commutative group scheme $G$
over~$K$.  We will now describe how to determine the structure of the
Abelian group $G(K)$.  After enlarging $K$ if necessary, we assume
that $K$ contains a root of unity $\zeta$ of order~$n$; this is not a
real restriction because one can take $\Gal(K(\zeta)/K)$-invariants
afterwards if desired.

Viewing the constant group scheme $G(K)_K$ over~$K$ as a subgroup
scheme of~$G$, we write $G(K)_K^\perp$ for the orthogonal complement
of $G(K)_K$ in $G^*$ under the duality pairing, and we put
$$
G^*\{K\} = G^*/G(K)_K^\perp,
$$
Then $G(K)_K$ corresponds to the largest quotient algebra $A'$ of~$A$
that is a product of copies of~$K$, and $G^*\{K\}$ corresponds to a
$K$-subalgebra $B'\subseteq B$.  We have now reduced to the case of
constant group schemes over~$K$, and we have a perfect pairing
$$
\pairing\colon G(K)\times G^*\{K\}\longrightarrow K^\times.
$$
We assume that we can explicitly determine a root of unity $\zeta\in
K$ of order~$n$ as well as all elements of
$G(K)\cong\Hom_{\Alg_K}(A,K)\simeq\Hom_{\Alg_K}(A',K)$ and
$G^*\{K\}\cong\Hom_{\Alg_K}(B',K)$.

\algorithm (Group structure). Given a field~$K$ containing a root of
unity $\zeta$ of order~$n$ and a dual pair $(A,B,\Phi)$ of
$K$-algebras of degree~$n$ over~$K$, this algorithm outputs the
sequence $d=(d_1,\ldots,d_r)$ of elementary divisors of $G(K)$ and
a group isomorphism
$$
f\colon H_d\isom\Hom_{\Alg_K}(A,K).
$$

\step Compute the largest quotient algebra $A'$ of~$A$ that is a
product of copies of~$K$.

\step Compute the subalgebra $B'\subseteq B$ as the orthogonal
complement of the kernel of the quotient map $A\to A'$ under~$\Phi$,
and $\Phi$ induces a perfect pairing $\Phi'\colon A'\times B'\to K$.

\step Fix bijections
$\alpha\colon\{1,\ldots,m\}\isom\Hom_{\Alg_K}(A',K)$ and
$\beta\colon\{1,\ldots,m\}\isom\Hom_{\Alg_K}(B',K)$.

\step Using the group isomorphism
$$
\lambda\colon\langle\zeta\rangle\isom{1\over n}\Z/\Z
$$
sending $\zeta$ to $1/n$, compute the matrix
$$
T=\bigl(\lambda(\langle\alpha(i),\beta(j)\rangle)\bigr)_{i,j=1}^m
$$
with entries in ${1\over n}\Z/\Z$.

\step By applying Algorithm~\ref{algorithm:group-structure} to~$T$,
determine a sequence $d$ of elementary divisors and bijections
$p\colon\{1,\ldots,m\}\isom H_d$ and $q\colon\{1,\ldots,m\}\isom
H_d^*$.

\step Output $d$ and the bijection
$$
\alpha\circ p^{-1}\colon H_d\isom\Hom_{\Alg_K}(A,L).
$$

\endalgorithm

\remark For this algorithm to work, it is in fact only necessary that
$K$ contains a root of unity of order equal to the exponent of $G(K)$.

\remark One can also use ``black box'' algorithms to determine the
group structure of $G(L)$; see for example Buchmann and
Schmidt \cite{Buchmann-Schmidt}.  However, the advantage of our
approach is that it only uses the pairing and does not need to perform
any group operations.

\subsection Validating input and determining group structures via
complex approximation

\label{subsec:complex}

We now describe a numerical variant of the algorithm
from \S\ref{subsec:group-structure} to verify whether a given triple
$(A,B,\Phi)$, where $A$ and $B$ are finite $\Q$-algebras with
distinguished $\Q$-bases and $\Phi$ is the matrix of a perfect
$\Q$-bilinear map $A\times B\to\Q$, is a dual pair of $\Q$-algebras.
In case the answer is yes, this algorithm also determines the group
structure of $G(\Qbar)$, where $G$ is the finite commutative group
scheme corresponding to $(A,B,\Phi)$.  Although our approach uses
numerical approximations, it is made rigorous thanks to height bounds.
It is important to stress that this approach does not require us to
compute any splitting fields.  We restrict to the field~$\Q$ and the
embedding $\Q\to\C$ for clarity; the idea can be generalised to other
number fields and to ultrametric places instead of Archimedean places.

If $K$ is a number field, we write $\Omega_K$ for the set of places
of~$K$; for every place $v$ of~$K$ we denote by $|\blank|_v\colon
K\to\R$ the normalised absolute value defined by~$v$.  Let
$h\colon\Qbar\to\R$ denote the absolute logarithmic height function,
given by
$$
h(x) = {1\over[K:\Q]}\sum_{v\in\Omega_K}\log\max\{1,|x|_v\}
$$
where $K$ is a number field with $x\in K\subset\Qbar$; this is
independent of the choice of~$K$.  It is straightforward to check that
for all $x,y\in\Qbar$ we have
$$
h(xy)\le h(x)+h(y)
\quad\hbox{and}\quad
h(x+y)\le h(x)+h(y)+\log 2.
$$

\proclaim Lemma. Let $\alpha\in\C$ be algebraic of degree at most~$d$.
If $\alpha\ne0$, then we have
$$
\bigl|-\log|\alpha|\bigr|\le dh(\alpha).
$$

\label{lemma:log-height}

\proof This follows from the identity
$$
h(\alpha)=h(1/\alpha),
$$
the inequality
$$
h(\alpha)\ge d^{-1}\log\max\{1,|\alpha|\}
\ge d^{-1}\log|\alpha|
$$
and the corresponding inequality for $1/\alpha$.\endproof

\proclaim Lemma. Let $\beta\in\C$ be algebraic of degree at most~$d$,
and let $\zeta\in\C$ be an $n$-th root of unity.  Then we have the
implication
$$
|\beta-\zeta| < \exp\bigl(-d\phi(n)(h(\beta)+\log 2)\bigr)
\;\Longrightarrow\; \beta=\zeta,
$$
where $\phi$ is Euler's $\phi$-function.

\label{lemma:close-to-zeta}

\proof This follows from Lemma~\ref{lemma:log-height} applied to
$\alpha=\beta-\zeta$; note that $\alpha$ lies in a number field of
degree at most $d\phi(n)$ and has height at most $h(\beta)+\log2$
since $h(\zeta)=0$.\endproof

To motivate the following algorithm, we note that if $(A,B,\Phi)$ is a
dual pair of $\Q$-algebras and $\sigma\colon A\to\C$ and~$\tau\colon
B\to\C$ are $\Q$-algebra homomorphisms, then $\sigma$ and~$\tau$
represent $\C$-valued points of the corresponding pair of Cartier dual
group schemes $G$ and~$G^*$.  Furthermore,
$(\sigma\otimes\tau)(\theta_\Phi)$ is the value of the duality pairing
between $\sigma$ and~$\tau$, and in particular is an $n$-th root of
unity.

\algorithm (Input validation and group structure using approximation).
Given two finite $\Q$-algebras $A$ and~$B$ of degree~$n$ with fixed
$\Q$-bases $(a_1,\ldots,a_n)$ and $(b_1,\ldots,b_n)$, respectively,
and the matrix of a perfect $\Q$-bilinear map $\Phi\colon A\times
B\to\Q$ with respect to these bases, this algorithm decides whether
$(A,B,\Phi)$ defines a commutative group scheme $G$ over~$\Q$.  If so,
it outputs the sequence $d=(d_1,\ldots,d_r)$ of elementary divisors of
$G(\Qbar)\cong G(\C)$ and group isomorphisms
$$
\eqalign{
f\colon H_d&\isom\Hom_{\Alg_\Q}(A,\C),\cr
g\colon H_d^*&\isom\Hom_{\Alg_\Q}(B,\C)}
$$
such that for all $x\in H_d$ and $\xi\in H_d^*$ we have
$$
\exp(2\pi i \xi(x)) = \langle f(x),g(\xi)\rangle,
$$
where the pairing on the right is induced by the canonical
isomorphisms $\Hom_{\Alg_\Q}(A,\C)\cong G(\C)$ and
$\Hom_{\Alg_\Q}(B,\C)\cong\Hom(G(\C),\C^\times)$.  Otherwise, it
outputs ``False''.

\step Check whether $A$ and~$B$ are \'etale over~$\Q$; if not, output
``False'' and stop.

\step Compute the matrix $\Theta$ of the element $\theta_\Phi\in
A\otimes_\Q B$ with respect to the fixed bases of $A$ and~$B$ as the
inverse transpose of the matrix of~$\Phi$.

\step Determine a positive integer $d$ such that the element
$\theta_\Phi\in A\otimes_\Q B$ is a root of a monic polynomial of
degree at most $d$ over~$\Q$ (for example, take $d=n^2$).

\step Fix a numbering $(\sigma_1,\ldots,\sigma_n)$ of the $\Q$-algebra
homomorphisms $A\to\C$ and a numbering $(\tau_1,\ldots,\tau_n)$ of the
$\Q$-algebra homomorphisms $B\to\C$.

\step Compute a numerical approximation to the $n\times n$-matrices
$$
P = (\sigma_i(a_j))_{i,j=1}^n,\quad
Q = (\tau_i(b_j))_{i,j=1}^n
$$
over~$\C$, with sufficient complex precision to perform the remaining
steps.

\step Compute a real number $C\ge 0$ such that for all $\Q$-algebra
homomorphisms $\psi\colon A\otimes_\Q B\to\Qbar$ we have
$h(\psi(\theta_\Phi))\le C$.

\step Compute a numerical approximation to the $n\times n$-matrix
$$
Z = P\Theta Q^{\rm t} = ((\sigma_i\otimes\tau_j)(\theta_\Phi))_{i,j=1}^n
$$
over~$\C$, with sufficient precision to perform the next step.

\step Use Lemma~\ref{lemma:close-to-zeta} to verify whether the
entries of~$Z$ are $n$-th roots of unity, and if so, to compute the
unique $n\times n$-matrix $T$ with coefficients in ${1\over n}\Z/\Z$
such that $Z$ is obtained from~$T$ by pointwise applying the map
$t\mapsto\exp(2\pi i t)$.  Otherwise, output ``False'' and stop.

\step Using Algorithm~\ref{algorithm:group-structure}, determine
whether $T$ describes an Abelian group.  If so, output $d$ and the
bijections $\alpha\circ p^{-1}$ and $\beta\circ q^{-1}$, where $d$,
$p$ and $q$ are the output of
Algorithm~\ref{algorithm:group-structure} and where
$\alpha\colon\{1,\ldots,n\}\isom\Hom_{\Alg_\Q}(A,\C)$ and
$\beta\colon\{1,\ldots,n\}\isom\Hom_{\Alg_\Q}(B,\C)$ are our fixed
numberings of the $\Q$-algebra homomorphisms $A\to\C$ and $B\to\C$.

\endalgorithm

\remark If one is willing to perform $\Q$-algebra operations in
$A\otimes_\Q B$, then one can first check whether $\theta_\Phi$ is an
$n$-th root of unity in $A\otimes_\Q B$, which is a necessary
condition for $(A,B,\Phi)$ to be a dual pair.  If so, one can take
$C=0$ and $d=\phi(n)$.  Moreover, instead of applying Lemma~5.2 to
find the matrix~$T$, it suffices to determine, for each entry of~$Z$,
the $n$-th root of unity closest to it.

\section Galois representations

\label{sec:galois-representations}

Let $K$ be a field, and let $\bar K$ be a separable closure of~$K$.
In this section, by a {\it Galois representation of~$K$\/} we will
mean a finite Abelian group $V$ of order not divisible by the
characteristic of~$K$ equipped with a continuous action of the Galois
group $\Gal(\bar K/K)$.  Given two Galois representations $V$ and~$W$
of~$K$, we view the Abelian group $\Hom(V,W)$ of all group
homomorphisms $V\to W$ as a Galois representation of~$K$ by equipping
it the with the $\Gal(\bar K/K)$-action defined by
$$
(\sigma f)(v) = \sigma(f(\sigma^{-1}v))
\quad\hbox{for all $v\in V$, $f\in\Hom(V,W)$ and $\sigma\in\Gal(\bar K/K)$}.
$$
We will apply this in particular to the case $W=\Gm(\bar K)$.

Given a Galois representation $V$ of~$K$, let $A_V$ denote the
finite \'etale $K$-algebra of $\Gal(\bar K/K)$-equi\-variant functions
$V\to\bar K$.  Then $A_V$ has a natural structure of Hopf algebra
over~$K$, and the finite \'etale $K$-scheme
$$
H_V = \Spec A_V
$$
has a natural structure of finite commutative group scheme over~$K$.
Conversely, one can reconstruct $V$ from either $A_V$ or $H_V$ by
$$
V = \Hom_{\Alg_K}(A_V,\bar K) = H_V(\bar K).
$$
We can therefore write down Galois representations in the form of dual
pairs of $K$-algebras.

\subsection Computing the matrix of an automorphism under a Galois
representation

Let $(A,B,\Phi)$ be a dual pair of $K$-algebras of degree~$n$, where
$n$ is not divisible by the characteristic of~$K$, and let $G$ be the
corresponding finite commutative group scheme over~$K$.  Let $L$ be an
extension field of~$K$ containing a root of unity $\zeta$ of
order~$n$.  In particular, if $L$ is a Galois extension of~$K$, then
$G(L)$ is a Galois representation of~$K$.  We will show how to
evaluate the map
$$
\Aut_K L\longrightarrow \Aut(G(L)).
$$
We assume that we have computed the structure of the group $G(L)$ as
in \S\ref{subsec:group-structure}.  In particular, we have a
sequence~$d=(d_1,\ldots,d_r)$ of elementary divisors and isomorphisms
$$
H_d\isom G(L),\quad
H_d^*\isom G^*\{L\}.
$$
As in \S\ref{subsec:group-structure}, let
$\lambda\colon\langle\zeta\rangle\isom{1\over n}\Z/\Z$ be the
isomorphism sending $\zeta$ to $1/n$.  Let $P_1,\ldots,P_r\in G(L)$
and $Q_1,\ldots,Q_r\in G^*\{L\}$ be the images of the standard
generators of~$H_d$ under these isomorphisms.  We also assume that the
matrix
$$
U = \bigl(\lambda(\langle P_i,Q_j\rangle)\bigr)_{i,j=1}^r
$$
with entries in ${1\over n}\Z/\Z$ is known.  (This matrix can in fact
easily be extracted from the data computed while determining the group
structure of $G(L)$.)

We will use the following notation.  If $d=(d_1,\ldots,d_r)$ is a
sequence of elementary divisors and $\End H_d$ is the (not necessarily
commutative) endomorphism ring of the Abelian group $H_d$, then we
identify the additive group of $\End H_d$ with the direct sum
$$
\bigoplus_{i=1}^r\bigoplus_{j=1}^r\Hom(\Z/d_j\Z,\Z/d_i\Z).
$$
Under this identification, an element $f\in \End H_d$ corresponds to a
collection of group homomorphisms $f_{i,j}\colon\Z/d_j\Z\to\Z/d_i\Z$.
We represent $f$ by the $r\times r$-matrix $(f_{i,j}(1))_{i,j=1}^r$;
note that $f_{i,j}(1)$ is an element of $\Z/d_i\Z$ annihilated
by~$d_j$.

Now let $\sigma$ be an automorphism of $L$ over~$K$.  We compute the
element $M(\sigma)\in\Aut H_d$ (represented by a matrix as above) as
follows.  We compute the matrix
$$
V(\sigma) = \bigl(\lambda(\langle\sigma P_i,Q_j\rangle)\bigr)_{i,j=1}^r
$$
with entries in ${1\over n}\Z/\Z$.  By the definition of $M(\sigma)$,
we have
$$
\sigma P_i = \sum_{k=1}^r M(\sigma)_{k,i} P_k.
$$
From the above identities, we deduce the system of equations
$$
V(\sigma)_{i,j} = \sum_{k=1}^r M(\sigma)_{k,i} U_{k,j},
$$
which we can solve for the $M(\sigma)_{k,i}$.

\example Assume that $K$ is a number field.  Let $\fp$ be a finite
place of~$K$ such that $(A,B,\Phi)$ has good reduction at~$\fp$, in
the sense that $(A,B,\Phi)$ can be written as the base change of a
dual pair $(A_\fp,B_\fp,\Phi_\fp)$ of algebras over the valuation ring
$\Z_{K,\fp}\subset K$ of~$\fp$.  Then we can compute the matrix of a
Frobenius element at~$\fp$ as follows.  Let $(\tilde A_\fp,\tilde
B_\fp,\tilde\Phi_\fp)$ be the base change of $(A_\fp,B_\fp,\Phi_\fp)$
to the residue field $k(\fp)$ of~$\fp$.  Applying the algorithm
described above to the dual pair $(\tilde A_\fp,\tilde
B_\fp,\tilde\Phi_\fp)$ of algebras over $k(\fp)$ and the Frobenius
automorphism of a finite splitting field of this dual pair
over~$k(\fp)$, we obtain the matrix of a Frobenius elements at~$\fp$
up to conjugacy in $\End H_d$.

\subsection Computing a dual pair corresponding to a Galois
representation

We now describe how to compute a dual pair of $K$-algebras
corresponding to a Cartier dual pair of Galois representations.
Suppose that we are in the situation where we have two Galois
representations $V$ and~$V'$ of~$K$ together with a perfect pairing
$$
\pairing\colon V\times V'\to\Gm(\bar K).
$$
that is $\Gal(\bar K/K)$-equivariant in the sense that
$$
\langle\sigma v,\sigma v'\rangle=\sigma\langle v,v'\rangle
\quad\hbox{for all $v\in V$, $v'\in V'$ and $\sigma\in\Gal(\bar K/K)$}
$$
Then the finite commutative group schemes $H_V$ and $H_{V'}$ over~$K$
are Cartier dual to each other.

We can compute a dual pair of $K$-algebras corresponding to these
group schemes as follows.  We choose $\Gal(\bar K/K)$-equivariant
functions
$$
\psi\colon V\to\bar K,\quad
\psi'\colon V'\to\bar K.
$$
Then the monic polynomials
$$
f=\prod_{v\in V}(x-\psi(v))\in\bar K[x],\quad
g=\prod_{v'\in V'}(y-\psi'(v'))\in\bar K[y]
$$
have coefficients in~$K$.  We write
$$
A=K[x]/(f),\quad B=K[y]/(g).
$$
Furthermore, we define an element
$$
\theta\in A\otimes_K B\cong K[x,y]/(f(x),g(y))
$$
as the unique element satisfying the interpolation property
$$
\theta(\psi(v),\psi'(v')) = \langle v,v'\rangle
\quad\hbox{for all $v\in V$, $v'\in V'$}.
$$
We can then compute the matrix of $\Phi$ with respect to the power
bases of $A$ and~$B$ as the inverse transpose of the matrix
of~$\theta$.

If $K$ is a number field and $V$ and~$V'$ are realised inside some
quasi-projective variety over~$K$, then except in the smallest
examples, the above computations are only feasible in practice using
numerical approximation.

\remark One can in fact obtain an explicit upper bound on the height
of the matrix of~$\Phi$ in terms of the heights of the roots of the
polynomials $f$ and~$g$.  We do not give the details here.

\subsection Explicit computations

The author explicitly computed several dozen dual pairs of
$\Q$-algebras corresponding to Galois representations of~$\Q$ on
groups of the form $\Z/n\Z$ or $\Z/n\Z\oplus\Z/n\Z$.  The sources of
these Galois representations are characters of $\Gal(\Qbar/\Q)$,
torsion subschemes of elliptic curves over~$\Q$, and Hecke eigenforms
over finite fields.

In practice, the data defining a dual pair $(A,B,\Phi)$ of
$\Q$-algebras often turns out to have smaller height than the data
defining the corresponding Hopf algebra, especially if $A$ and~$B$ are
not isomorphic (i.e.\ the group scheme is not isomorphic to its own
Cartier dual).  This is useful in practice because our computations
rely on numerical approximation, and the required precision is
directly related to the height of the output.

\example The author computed a dual pair of $\Q$-algebras
corresponding to the Cartier dual pair of Galois representations
$$
\rho_1,\rho_2\colon\Gal(\Qbar/\Q)\to\GL_2(\F_7)
$$
attached to the two Hecke eigenforms $f_1$ and~$f_2$ for the
group~$\Gamma_1(13)$ over~$\F_7$ with $q$-expansions
$$
\eqalign{
f_1 &= q + q^2 + q^3 + 5q^4 + 5q^5 + q^6 + 2q^8 + 2q^9 + O(q^{10}),\cr
f_2 &= q + 3q^2 + 4q^3 + 3q^4 + 2q^5 + 5q^6 + 5q^8 + 4q^9 + O(q^{10}).}
$$
These representations occur in the 7-torsion of the Jacobian of the
modular curve $\X_1(13)$ of genus~2; the duality between them is
realised by the Weil pairing.  The computations were done using the
author's {\tt modgalrep} software package \cite{modgalrep}, written in
C and PARI/GP \cite{pari}.  The height of the resulting polynomials
$f\in\Q[x]$ and $g\in\Q[y]$ (i.e.\ the maximum of the logarithmic
heights of the coefficients) is about 88 for~$f$ and 93 for~$g$.  The
height of the matrix of~$\Phi$ with respect to the power bases of
$A=\Q[x]/(f)$ and~$B=\Q[y]/(g)$ is about 585.  In contrast, the height
of the image of~$x$ under the comultiplication map
$$
\mu\colon A\longrightarrow A\otimes_\Q A
\cong\Q[x_1,x_2]/(f(x_1),f(x_2))
$$
is roughly 1626.

\section Future work

\label{sec:future}

As we have seen, dual pairs of algebras over a ring~$R$ are a concise
way to write down finite locally commutative group schemes over~$R$,
and the category of dual pairs of algebras over~$R$ is equivalent to
the category of finite commutative group schemes over~$R$.  Moreover,
if $R$ is a number field, the resulting objects have comparatively
small height.  For these reasons, in the near future the author is
planning to work on building a database of Galois representations
stored as dual pairs of algebras.

Besides representing finite commutative group schemes by a dual pair
of algebras over a ring~$R$, it is also possible to represent {\it
torsors\/} over such group schemes.  If $G$ is a finite commutative
group scheme represented by the dual pair $(A,B,\Phi)$ over~$R$, then
as an algebraic structure representing $G$-torsors one can take a
triple $(T,U,\Psi)$, where $T$ is a finite locally free $R$-algebra,
$U$ is a locally free $B$-module of rank~$1$, and $\Psi\colon T\times
U\to R$ is a perfect $R$-bilinear map.  Another future direction is to
develop algorithms for working with torsors in this representation;
this could be useful for computing in Selmer groups.

\section Appendix: identifying a finite Abelian group from a pairing

Let $H$ be a finite Abelian group of order~$n$, and let
$H^*=\Hom(H,\Q/\Z)$ be the dual group.  Let
$p\colon\{1,2,\ldots,n\}\isom H$ and $q\colon\{1,2,\ldots,n\}\isom
H^*$ be two bijections, viewed as enumerations of the elements of $H$
and~$H^*$, respectively.  To these data we associate the $n\times
n$-matrix $T=(q(j)(p(i)))_{i,j=1}^n$ with entries in the group
${1\over n}\Z/\Z$.

We say that an $n\times n$-matrix $T$ with entries in ${1\over
n}\Z/\Z$ {\it describes an Abelian group\/} if it is of the form
$(q(j)(p(i)))_{i,j=1}^n$ for some choice of $H$, $p$ and~$q$ as above.
It is not hard to see that $T$ determines $H$ up to isomorphism.

Algorithm~\ref{algorithm:group-structure} below checks whether a given
$n\times n$-matrix with entries in ${1\over n}\Z/\Z$ describes an
Abelian group~$H$, and if so, identifies $H$.  For the latter, we will
use the following terminology and notation.  By a {\it sequence of
elementary divisors\/} we mean a sequence $d=(d_1,\ldots,d_r)$ of
integers greater than~$1$ satisfying $d_r\mid d_{r-1}\mid\ldots\mid
d_1$.  For such a sequence~$d$, we write
$$
H_d = \Z/d_1\Z\oplus\cdots\oplus\Z/d_r\Z.
$$

\algorithm (Identify a finite Abelian group from a pairing matrix).
Given a positive integer~$n$ and an $n\times n$-matrix $T$ with
entries in ${1\over n}\Z/\Z$, this algorithm decides whether $T$
describes an Abelian group.  If so, it outputs a sequence~$d$ of
elementary divisors and two bijections
$$
p\colon\{1,2,\ldots,n\}\isom H_d
\quad\hbox{and}\quad
q\colon\{1,2,\ldots,n\}\isom H_d^*
$$
satisfying $T=(q(j)(p(i)))_{i,j=1}^n$; otherwise, it outputs
``False''.

\step If $n=1$, output the empty sequence $d$ and the unique
bijections $p$ and~$q$.

\step Let $d_1$ be the maximum of the denominators of the entries
of~$T$.  If $T$ contains the element $1/d_1$, then choose a pair
$(i_1,j_1)$ of row and column indices such that $T_{i_1,j_1}=1/d_1$;
otherwise, output ``False'' and stop.

\step Determine injective functions
$$
f_1\colon\Z/d_1\Z\to\{1,2,\ldots,n\},\quad
g_1\colon\Z/d_1\Z\to\{1,2,\ldots,n\}
$$
such that for all $x\in\Z/d_1\Z$, the $f_1(x)$-th row of~$T$ equals
$x$ times the $i_1$-th row and the $g_1(x)$-th column of~$T$ equals
$x$ times the $j_1$-th column.  If these do not exist, output
``False'' and stop.

\step Write $n'=n/d_1$.  Let $T'$ be the submatrix formed by the
$T_{i,j}$ where $i$ runs over the row indices such that $T_{i,j_1}=0$
and $j$ runs over the column indices such that $T_{i_1,j}=0$.  If $T'$
is not an $n'\times n'$-matrix with entries in ${1\over n'}\Z/\Z$,
output ``False'' and stop.

\step Apply the algorithm recursively to $T'$ to check whether $T'$
describes an Abelian group, and if so, to find the corresponding
sequence $d'$ of elementary divisors and bijections $p'$ and~$q'$.
Since $T'$ is a submatrix of~$T$, the inverses of $p'$ and~$q'$ can be
viewed as injective functions
$$
f'\colon H_{d'}\to\{1,2,\ldots,n\},\quad
g'\colon H_{d'}^*\to\{1,2,\ldots,n\}.
$$

\step Writing $d'=(d_2,\ldots,d_r)$, let $d=(d_1,d_2,\ldots,d_r)$.
Determine whether $d_2$ divides $d_1$, and whether there exist
(necessarily unique) bijections
$$
f\colon H_d\isom\{1,2,\ldots,n\},\quad
g\colon H_d^*\isom\{1,2,\ldots,n\}
$$
such that for all $x_1\in\Z/d_1\Z$ and $x'\in H_{d'}$, the
$f(x_1,x')$-th row of~$T$ equals the $f_1(x_1)$-th row plus the
$f'(x')$-th row, and similarly for $g$ (using the columns of~$T$).  If
so, put $p=f^{-1}$ and $q=g^{-1}$ and output $(d,p,q)$; otherwise,
output ``False'' and stop.

\endalgorithm

\label{algorithm:group-structure}

\unnumberedsection References

\normalbaselines
\rm
\normalparindent=25pt
\parindent=\normalparindent
\parskip=1ex plus 0.5ex minus 0.2ex

\reference{Bosman} J. G. {\sc Bosman}, {\sl Explicit computations with
modular Galois representations\/}.  Ph.\thinspace D.\ thesis,
Universiteit Leiden, 2008.

\reference{thesis} P. J. {\sc Bruin}, {\sl Modular curves, Arakelov
theory, algorithmic applications\/}.  Ph.\thinspace D.\ thesis,
Universiteit Leiden, 2010.

\reference{modgalrep} P. J. {\sc Bruin}, {\tt modgalrep}, software
package for computing Galois representations attached to modular forms
over finite fields.  Available at {\tt
http://gitlab.com/pbruin/modgalrep/}.

\reference{Buchmann-Schmidt} J. {\sc Buchmann} and A. {\sc Schmidt}, 
Computing the structure of a finite abelian group.  {\it Mathematics
of Computation\/} {\bf 74} (2005), no.~252, 2017--2026.

\reference{DvHZ} M. {\sc Derickx}, M. {\sc van Hoeij} and J. {\sc
Zeng}, Computing Galois representations and equations for modular
curves $X_H(\ell)$.  Preprint, {\tt https://arxiv.org/abs/1312.6819}.

\reference{Edixhoven-Couveignes} {\sl Computational aspects of modular
forms and Galois representations\/}.  Edited by S. J. {\sc Edixhoven}
and J.-M. {\sc Couveignes}, with contributions by J. G. {\sc Bosman},
J.-M. {\sc Couveignes}, S. J. {\sc Edixhoven}, R. S. {\sc de Jong}
and F. {\sc Merkl}.  Annals of Mathematics Studies {\bf 176}.
Princeton University Press, Princeton, NJ, 2011.

\reference{Hartshorne} R. {\sc Hartshorne}, {\sl Algebraic
Geometry\/}.  Springer-Verlag, New York, 1977.

\reference{Mac Lane} S. {\sc Mac Lane}, {\sl Categories for the Working
Mathematician\/}.  Second edition.  Graduate Texts in Mathematics {\bf
5}.  Springer-Verlag, New York, 1998.

\reference{Mascot} N. {\sc Mascot}, Computing modular Galois
representations.  {\it Rendiconti del Circolo Matematico di Palermo
(2)\/} {\bf 62} (2013), no.~3, 451--476.

\reference{Oort} F. {\sc Oort}, {\sl Commutative group schemes\/}.
Lecture Notes in Mathematics {\bf 15}.  Springer-Verlag,
Berlin/Heidelberg/New York, 1966.

\reference{pari} The PARI~Group, PARI/GP, version {\tt 2.9.0}.
Universit\'e de Bordeaux, 2016,\hfill\break
{\tt http://pari.math.u-bordeaux.fr/}.

\reference{Quillen} D. {\sc Quillen}, Higher algebraic $K$-theory: I.
In: H. {\sc Bass} (editor), {\sl Higher $K$-Theories\/}.  Lecture
Notes in Mathematics {\bf 341}.  Springer-Verlag,
Berlin/Heidelberg/New York, 1973.

\reference{Tate} J. {\sc Tate}, The arithmetic of elliptic curves.
{\it Inventiones mathematicae\/} {\bf 23} (1974), 179--206.

\reference{Tate-Oort} J. {\sc Tate} and F. {\sc Oort}, Group schemes
of prime order.  Annales scientifiques de l'\'E.N.S. ($4^{\rm e}$
s\'erie) {\bf 3} (1970), 1--21.

\reference{Tian} P. {\sc Tian}, Computations of Galois representations
associated to modular forms of level one.  {\it Acta Arithmetica\/}
{\bf 164} (2014), no.~4, 399--412.

\reference{Yin-Zeng} L. {\sc Yin} and J. {\sc Zeng}, On the
computation of coefficients of modular forms: the reduction modulo~$p$
approach.  {\it Mathematics of Computation\/} {\bf 84} (2015),
no.~293, 1469--1488.

\vskip 2cm

\leftline{Peter Bruin}
\leftline{Universiteit Leiden}
\leftline{Mathematisch Instituut}
\leftline{Postbus 9512}
\leftline{2300 RA \ Leiden}
\leftline{Netherlands}
\leftline{\tt P.J.Bruin@math.leidenuniv.nl}

\bye